\documentclass[titlepage,11pt]{article}
\usepackage{latexsym}
\usepackage{amsfonts}
\usepackage{amsmath}
\usepackage{mathtools}
\usepackage{comment}
\usepackage{tikz}
\usetikzlibrary {shapes.geometric}
\newcommand\blackslug{\hbox{\hskip 1pt \vrule width 4pt height 8pt depth 1.5pt
        \hskip 1pt}}
\newcommand\bbox{\hfill \quad \blackslug \bigbreak}

\def\LL{,\ldots,}

\newcommand{\cupcup}{\cup \cdots\cup}

\newcommand{\mac}{\mathcal}

\def\dist{\operatorname{dist}}
\def\bd{\operatorname{bd}}
\newcommand{\bS}{\bd(\Sigma)}

\title{Asymptotic structure. VI. Distant paths across across a disc}
\author{
Tung Nguyen\thanks{Supported by a Porter Ogden Jacobus Fellowship. Current address: University of Oxford,
Oxford, UK.}\\
Princeton University,\\ Princeton, NJ 08544, USA
\and
Alex Scott\thanks{Supported by EPSRC grant EP/X013642/1}\\
University of Oxford, \\
Oxford, UK
\and
Paul Seymour\thanks{Supported by AFOSR grant
FA9550-22-1-0234, and by NSF grant  DMS-2154169.}\\
Princeton University,\\ Princeton, NJ 08544, USA}

\date{July 26, 2025; revised \today}

\newtheorem{thm}{}[section]

\newcommand{\Proof}{\noindent{\bf Proof.}\ \ }

\begin{document}
\maketitle
\begin{abstract}
Menger's theorem says that, for $k\ge0$, if $S,T$ are sets of vertices in a graph $G$, then either there are $k+1$ vertex-disjoint paths between $S$ and $T$, or there is a set  $X$ of at most $k$ vertices such that every $S$-$T$ path passes through $X$.  The  ``coarse Menger conjecture'' proposed a generalization of Menger's theorem for paths that are far apart: for all $k,c$ there exists $\ell$, such that for every graph $G$ and subsets $S,T\subseteq V(G)$, either there are $k+1$ paths between $S$ and $T$, pairwise with distance more than $c$, or there is a set $X\subseteq V(G)$ of at most $k$ vertices such that every $S$-$T$ path has distance at most $\ell$ from $X$. 
This is known to be false, but may be true if $G$ is planar. Here we show that it is true if $G$ is planar and all vertices in 
$S\cup T$ are on the infinite region.  In this case, we also 
obtain a linear-time algorithm to test for the existence of $k + 1$ paths between $S$ and $T$, pairwise with distance more than $c$.
\end{abstract}

\section{Introduction}

The ``disjoint paths problem'' asks when there is a set of $k+1$ vertex-disjoint paths between sets $S,T$ of vertices of a graph $G$; 
and it is answered by a theorem of Menger from 1927~\cite{menger}, that such paths exist if and only if there is no subset $X\subseteq V(G)$ 
of size at most $k$ such that every $S$-$T$ path has a vertex in $X$. But what if we want the paths to be at least a certain distance from one 
another?
This question is motivated both by the developing area of ``coarse graph theory'', which is concerned with the large-scale geometric 
structure of graphs (see Georgakopoulos and Papasoglu~\cite{agelos}), and by the algorithmic question of deciding whether such paths exist (see
Bienstock~\cite{bienstock}, Kawarabayashi and Kobayashi~\cite{kk}, and Balig\'acs and MacManus~\cite{bm}).  The purpose of this paper is to begin the investigation of whether a coarse analogue of Menger's theorem holds for graphs embeddable on a fixed surface. 

All graphs in this 
paper are finite, and have no loops or parallel edges.
If $X$ and $Y$ are vertices, or sets of vertices, or subgraphs, of a graph $G$, then $\dist_G(X,Y)$ denotes the distance between 
$X,Y$, that is, the number of edges in the shortest path of $G$ with one end in $X$ and the other in $Y$.
We say $X,Y$ are {\em $c$-distant} if their distance is at least $c$.

A coarse analogue of Menger's theorem was conjectured by
Albrechtsen, Huynh, Jacobs, Knappe and Wollan~\cite{wollan}, and independently by Georgakopoulos and Papasoglu~\cite{agelos}:
\begin{thm}\label{conj}
{\bf False conjecture (Coarse Menger conjecture):} For all integers $k,c\ge 0$ there exists $\ell>0$ with the following property.
Let $G$ be a graph and let $S,T\subseteq V(G)$.  Then either 
\begin{itemize}
\item there are $k+1$ paths between $S,T$, pairwise $c$-distant; or 
\item there is a set $X\subseteq V(G)$
with $|X|\le k$ such that every path between $S,T$ contains a vertex with distance at most $\ell$ from some member of $X$.
\end{itemize}
\end{thm}
Both sets of authors \cite{agelos,wollan} proved the conjecture in the case $k=1$; 
Gartland, Korhonen and Lokshtanov \cite{gkl} and Hendrey, Norin, Steiner, and Turcotte \cite{hnst} proved the case $c=2$ for graphs with bounded degree;
and we proved recently that the conjecture holds for graphs of bounded path-width \cite{as5}.
However, we showed in~\cite{counterex} that the coarse Menger conjecture is false for all $k\ge 2$,
even if $c=3$ and $G$ has maximum degree three.  Indeed, even a weak version of the conjecture, where we allow the size of $X$ to be 
an arbitrary function of $k$, turns out to be false~\cite{AS4}.

Our counterexamples are highly non-planar (they have unbounded genus), and one might hope that the Coarse Menger conjecture is true
for planar graphs, or perhaps even for graphs drawn in any fixed surface. In other words, perhaps the following is true.

\begin{thm}\label{newconj}
{\bf (Coarse Menger conjecture for surfaces):} Let $\mathcal S$ be a fixed surface. For all integers $k,c\ge 0$ there exists $\ell>0$ with the following property.
Let $G$ be a graph that embeds on $\mathcal S$ and let $S,T\subseteq V(G)$.  Then either 
\begin{itemize}
\item there are $k+1$ paths between $S,T$, pairwise $c$-distant; or 
\item there is a set $X\subseteq V(G)$
with $|X|\le k$ such that every path between $S,T$ contains a vertex with distance at most $\ell$ from some member of $X$.
\end{itemize}
\end{thm}

In this paper we begin the investigation of this question for planar graphs.  We show that it is true if $G$ is planar 
and can be drawn in the plane such that all vertices in $S\cup T$ are incident with the infinite region (briefly, ``with $S,T$ on the outside''). 
More exactly:
\begin{thm}\label{mainthm}
Let $c,k\ge 0$, and let $G$ be drawn in the plane, with $S,T$ on the outside.
Then either:
\begin{itemize}
\item there are $k+1$ paths between $S,T$, pairwise $(c+1)$-distant; or
\item there is a set of at most $k$ connected subgraphs of $G$,
such that every path between $S,T$ intersects one of these subgraphs, and the sum of the diameters of the subgraphs is at most
$200k^3c$.
\end{itemize}
\end{thm}
(Note that we have switched from ``$c$-distant'' to ``$(c+1)$-distant'', which is slightly more convenient.)
To prove this, we first obtain a necessary and sufficient condition
for the existence of $k$ pairwise $(c+1)$-distant paths $P_1\LL P_k$, where $P_i$ joins two given vertices $s_i,t_i$,
and the vertices $s_1\LL s_k,t_1\LL t_k$ all belong to the infinite region.

Finally, 
we obtain a linear-time algorithm to test whether there is a set of $k+1$ pairwise $(c+1)$-distant paths between $S$ and $T$ 
(for fixed $c,k$, and still assuming that $G$ is drawn in the plane and $S, T$ are on the outside of the drawing).

\section{The linkage problem}\label{sec:linkage}

The linkage problem is different from the question answered by Menger's theorem: now we are given $k$ pairs $(s_i,t_i)\;(1\le i\le k)$
of vertices of a graph $G$, and a {\em linkage} for this set of pairs means a set of $k$ vertex-disjoint paths $P_1\LL P_k$, where $P_i$ joins $s_i,t_i$ for 
$1\le i\le k$. Testing whether a linkage exists is NP-complete if $k$ is part of the input~\cite{karp}, and solvable in polynomial 
time if $k$ is fixed~\cite{GM13}. But it is much simpler if $G$ is drawn in a closed disc $\Sigma$ and $s_1\LL s_k, t_1\LL t_k$
belong to the boundary $\bd(\Sigma)$: it is shown in~\cite{GM6} that in that case, a linkage for the pairing $(s_i,t_i)\;(1\le i\le k)$ exists 
if and only if:
\begin{itemize}
\item  no two of the pairs $(s_i,t_i)$ 
``cross'' (that is, for $1\le i<j\le k$, there is a line segment of $\bd(\Sigma)$ that contains both of $s_i,t_i$ and neither of $s_j,t_j$); and
\item for every two points $a,b\in \bd(\Sigma)$ and every simple curve $L$ between $a,b$ in $\Sigma$ that meets the drawing only in vertices, 
$|L\cap V(G)|$  is at least the number of values of $i\in \{1\LL k\}$ 
such that both line segments of $\bd(\Sigma)$
between $a,b$ contain one of $s_i,t_i$.
\end{itemize}

In this paper we are asking for  a linkage of $(c+1)$-distant paths, but it turns out that, 
in the same disc case, a similar theorem holds; and this is a key lemma for the proof of \ref{mainthm}. In this section we state and prove it. 

If $G$ is drawn in $\Sigma$, each vertex of $G$ is a point of $\Sigma$ and each edge of $G$ is a line segment, that is, 
a subset of $\Sigma$
homeomorphic to the closed interval $[0,1]$. We denote the union of the set of vertices and the set of
edges of $G$ by $U(G)$. A {\em region} of $G$ in $\Sigma$ means an arc-wise connected component of $\Sigma\setminus U(G)$.

Now us fix some number $c\ge 0$ (throughout the paper); and again, let $G$ be drawn in a closed disc $\Sigma$.
A {\em log} is a path of $G$ of length at most $c$. A {\em boom}\footnote{Historically, a boom was an obstruction, made by logs chained together and strung across a harbour mouth, so it seems an appropriate name.} of {\em length $t$}
is a sequence $Q_1\LL Q_t$ of logs, such that for $1\le i<t$, there is a region of $G$
incident with a vertex of $Q_i$ and with a vertex of $Q_{i+1}$ (this is true, for instance, if $Q_i\cap Q_{i+1}$ is non-null).
If $b\in \bd(\Sigma)$, a boom $Q_1\LL Q_t$ {\em attaches to } $b$ if either:
\begin{itemize}
\item $b\in V(Q_1\cupcup Q_t)$, or
\item $b$ belongs to a region of $G$ that is incident with a vertex of $Q_1\cupcup Q_t$.
\end{itemize}
Note that the regions are disjoint from $U(G)$, so if $b$ satisfies the first bullet it is a vertex, and if it satisfies the second then
it belongs to $\bS \setminus U(G)$.

\begin{figure}[h!]
\centering

\begin{tikzpicture}[scale=.5,auto=left]

\tikzstyle{every node}=[inner sep=1.5pt, fill=black,circle,draw]
\draw[dotted, thick] (0,0) circle (5);

\node (v1) at (0,4) {};
\node  (v2) at (0,3.5) {};
\node (v3) at (0,3) {};
\node (v4) at (0,2.5) {};
\node (w1) at (.2,1) {};
\node (w2) at (.3, .5) {};
\node (w3) at (.4, 0) {};
\node (w4) at (.5, -0.5) {};
\node (x1) at (.3,-3) {};
\node (x2) at (.1,-3.5){};
\node (x3) at (-.1,-4) {};
\node (x4) at (-.3,-4.5){};
\node (x5) at (-.5,-5) {};

\draw (v1) -- (v4);
\draw (w1) -- (w4);
\draw (x1) -- (x5);

\draw [dotted, thick, fill = gray!50] (v4) to [bend right = 60] (w1) to [bend right = 60] (v4);
\draw [dotted, thick, fill = gray!50] (w4) to [bend right = 60] (x1) to [bend right = 60] (w4);
\draw [dotted, thick, fill = gray!50] (-.5, 5) to  [bend right = 40] (v1) to [bend right = 40] (.5, 5) to [bend right = 5] (-.5, 5);
\node (v1) at (0,4) {};
\node (v4) at (0,2.5) {};
\node (w1) at (.2,1) {};
\node (w4) at (.5, -0.5) {};
\node (x1) at (.3,-3) {};

\node (s1) at ({5*cos(150)}, {5*sin(150)}) {};
\node (s2) at ({5*cos(170)}, {5*sin(170)}) {};
\node (s3) at ({5*cos(190)}, {5*sin(190)}) {};
\node (s4) at ({5*cos(210)}, {5*sin(210)}) {};
\node (t1) at ({5*cos(30)}, {5*sin(30)}) {};
\node (t2) at ({5*cos(10)}, {5*sin(10)}) {};
\node (t3) at ({5*cos(-10)}, {5*sin(-10)}) {};
\node (t4) at ({5*cos(-30)}, {5*sin(-30)}) {};

\tikzstyle{every node}=[]
\draw[left] (s1) node []           {$s_1$};
\draw[left] (s2) node []           {$s_2$};
\draw[left] (s3) node []           {$t_3$};
\draw[left] (s4) node []           {$t_4$};
\draw[right] (t1) node []           {$t_1$};
\draw[right] (t2) node []           {$t_2$};
\draw[right] (t3) node []           {$s_3$};
\draw[right] (t4) node []           {$s_4$};

\end{tikzpicture}

\caption{A boom of length three, obstructing a linkage for $(s_i,t_i)\;(1\le i\le 4)$ of pairwise 4-distant paths. The gray areas are regions.}
\label{fig:boom}
\end{figure}
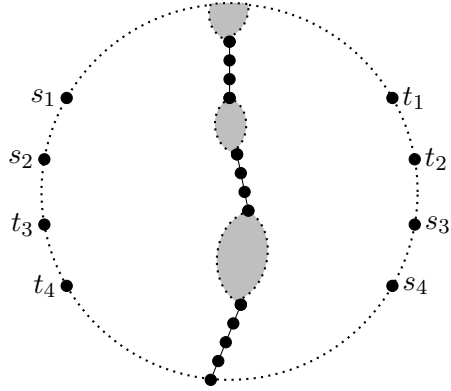

If $a,b\in \bd(\Sigma)$, a boom $Q_1\LL Q_t$ {\em joins} $a,b$ if either it attaches to both $a,b$, or $t=0$ and $a,b$ belong to a common region.
We will prove the following (see Figure \ref{fig:boom}):

\begin{thm}\label{linkage}
Let $\Sigma$ be a closed disc, let $G$ be a graph drawn in $\Sigma$, and let $c\ge 0$ be an integer.
Let $s_1,t_1\LL s_k, t_k\in V(G)\cap \bS$, such that $(s_i,t_i),(s_j,t_j)$ do not cross for $1\le i<j\le k$.
Then exactly one of the
following holds:
\begin{itemize}
\item There is a linkage for the pairs $(s_1,t_1)\LL (s_k,t_k)$, of paths that are pairwise $(c+1)$-distant;
\item there exist distinct $a,b\in \bd(\Sigma)$, and a boom joining $a,b$ of length less than
the number of $i\in \{1\LL k\}$ such that
$(s_i,t_i)$ crosses $(a,b)$.
\end{itemize}
\end{thm}

There is a problem that is slightly more general than the linkage problem across a disc described above, and it is useful
for our application to prove a result for this more complicated version. Let $C$ be a circle. A subset of $C$ that is nonempty, closed and arcwise connected is called an {\em interval} of $C$; so intervals are either line segments, or single points or the whole of $C$. Two intervals are {\em internally-disjoint} if no internal 
point of either of them belongs to the other. If $I_1,I_2,I_3,I_4$ are pairwise internally-disjoint intervals of $C$, we say 
that $(I_1,I_3)$
{\em crosses} $(I_2,I_4)$ if  there is a rotation around $C$ that traverses the sets $I_1,I_2,I_3,I_4$ in this order. (This extends 
the definition for crossing of pairs of points given earlier.) For intervals $I$ with $|I|=1$ and $|I|=\{x\}$ say, we often 
write $x$ for $I$.

Now let $\Sigma$ be a closed disc, and let us assign a direction of rotation ``clockwise'' to $\bd(\Sigma)$.  
If $A,B$ are distinct, internally-disjoint intervals of $\bS$, $[A\to B]$ denotes the interval $L$ of $\bS$ that is minimal subject to 
$L\cap A, L\cap B\ne \emptyset$, and such that $A,L,B$ are in clockwise order. 
Let $\mac I =\{I_1\LL I_k\}$ be a set of pairwise internally-disjoint intervals of $\bd(\Sigma)$, numbered in clockwise order. 
We call $\mac I=\{I_1\LL I_k\}$ an {\em interval system} for $\bd(\Sigma)$.
For $1\le i<j\le k$, let $d_{i,j}\ge 0$ be an integer, such that if $1\le p,q,r,s\le k$ then one of $d_{p,r}, d_{q,s}=0$.
We call the function $d$ a {\em demand function} for $\mac I$. 
If $G$ is drawn in $\Sigma$, a {\em linkage in $G$ 
for  $d,\mac I$} is a family $\mac L=(\mathcal{L}_{i,j}:1\le i<j\le k)$ of pairwise disjoint sets of paths, 
where $\mathcal{L}_{i,j}$
contains $d_{i,j}$ paths, all between $I_i, I_j$, and the paths in $U(\mac L)$ are pairwise vertex-disjoint
where $U(\mac L)$ denotes
the set of paths $\bigcup_{1\le i<j\le k}\mac L_{ij}$. (The condition ``one of $d_{p,r}, d_{q,s}=0$'' is assumed because there could not be a corresponding linkage if this condition was not satisfied.)
A linkage $\mac L$ is {\em $c$-distant} if every two members of $U(\mac L)$ are $c$-distant.

The next result implies \ref{linkage} by taking each $I_i$ to be a singleton.
We will use this, rather than \ref{linkage}, in the proof of \ref{mainthm}. 
\begin{thm}\label{intlinkage}
Let $\Sigma$ be a closed disc, let $\mac I = \{I_1\LL I_k\}$ be an interval system, and  let $d$ be a demand function for $\mac I$.
Let $G$ be a graph drawn in $\Sigma$, and let $c\ge 0$ be an integer.
Then exactly one of the 
following holds:
\begin{itemize}
\item there is a $(c+1)$-distant linkage in $G$ for $d, \mac I$;
\item there exist $a,b\in \bd(\Sigma)$, and a boom joining $a,b$, of length 
less than the sum of $d_{i,j}$ over all $1\le i<j\le k$ such that 
$(I_i,I_j)$ crosses $(a,b)$.
\end{itemize}
\end{thm}
\Proof It is easy to see that not both hold and we omit the argument. Now we prove that one of the statements holds.
We suppose for a contradiction that neither holds.
By slightly enlarging $\Sigma$, we may assume that no edge of $G$ intersects $\Sigma$, and the only vertices of $G$ in $\Sigma$
are the vertices in $I_1\cupcup I_k$. (This is not quite trivial: the enlargement does not change whether the first outcome of the 
theorem holds, but in the second outcome one of $a,b$, say $x$, might belong to $V(G)\setminus (S\cup T)$ and to the boundary of the old disc,  and needs to be replaced with 
some point in a region with $x$ of the boundary of the new disc.) We may also assume that $k\ge 2$, and so none of the intervals in $\mac I$ is the whole of $\bS$.
Moreover, we can assume that  for each of the intervals $I_i$, both its ends (or its point, if it is a singleton) belong to $V(G)$. 

If $P,P'$ are paths of $G$, with ends $s,t$ and $s',t'$ respectively, and $s,s',t,t'\in \bS$, we say that $P'$ is {\em inside} $P$ if
$P,P'$ are vertex-disjoint and  either 
\begin{itemize}
\item $s\ne t$ and $s',t'\in [s\to t]$; or
\item $s=t$ and $[I_k\to I_1]= V(P)$.
\end{itemize}

Now let $0\le d_{ij}'\le d_{ij}$ for $1\le i<j\le k$, such that for $1\le i\le j\le k$, if $d'_{ij}>0$ then $d'_{i'j'}=d_{i'j'}$
for all $i',j'$ with $i\le i'<j'\le j$ and $(i',j')\ne (i,j)$. We call $d'$ a {\em partial $d$-demand}. 
Let $\mac{L}=(\mathcal{L}_{ij}:1\le i<j\le k)$ be a $(c+1)$-distant linkage 
for the demand function $d'$ for $\mac I$. For $1\le i<j\le k$ and $P\in L_{ij}$, 
we say that $P$ is {\em pushed}
if 
for each vertex $v\in V(P)$, either:
\begin{itemize}
\item $v\in [I_i\to I_j]$; or
\item $v$ is incident with a region that contains a point of $[I_i\to I_j]$; or
\item there is a region incident with $v$ and some vertex $u$, and a log $Q$ 
with ends $u,w$, where $w$ belongs to 
some path $P'\ne P$ in $U(\mac L)$ that is inside $P$.
\end{itemize}
We say that $\mac L$ is 
{\em pushed} if for 
$1\le i<j\le k$, each $P\in \mac L_{i,j}$ is pushed.

Let us choose 
a partial $d$-demand $d'$, with $\sum d'_{ij}$ maximum such that there is a pushed $(c+1)$-distant linkage 
$\mac L=(\mac L_{ij}:1\le i<j\le k)$ for $d'$. We may further assume that for $1\le i<j\le k$ and each $P\in \mac L_{ij}$, $V(P)\cap I_i$
contains only one vertex, an end of $P$, and the same for $V(P)\cap I_j$. 
For every two internally-disjoint  intervals $A,B$ of $\bd(\Sigma)$, let $f(A,B)$ be the number of paths in $U(\mac L)$ such that its pairs of ends crosses $(A,B)$. 
\\
\\
(1) {\em For $1\le i<j\le k$, and each $P\in \mac L_{i,j}$, 
and for each $v\in V(P)$, there exists $b$ in $[I_i\to I_j]$ and not in the interior of any member of $\mac I$, 
and a boom $Q_1\LL Q_t$
attaching to $b$ and with $V(Q_1)=\{v\}$, of length $t=f(V(P)\cap I_i,b)$.}
\\
\\
Let $P$ have ends $s\in I_i$ and $t\in I_j$. 
Define $\phi(P)$ by:
\begin{itemize}
\item if $s\ne t$, let $\phi(P)$ be the number of paths in $U(\mac L)$ with both ends in $[s\to t]$ (including $P$ itself); 
\item if $s=t$ and $V(P)=[I_i \to I_j]$, let $\phi(P) = 1$; and
\item if $s=t$ and $V(P)=[I_j \to I_i]$, let $\phi(P) = |U(\mac L)|$.
\end{itemize}
We proceed by induction on $\phi(P)$. 
Let $v\in V(P)$. Since $P$ is pushed,  either:
\begin{itemize}
\item $v\in [I_i\to I_j]$; or
\item $v$ is incident with a region that contains a point of $[I_i\to I_j]$; or
\item there is a region incident with $v$ and some vertex $u$, and a log $Q$ with ends $u,w$, where $w$ belongs to some path $P'$ of 
$\mac L$
inside $P$.
\end{itemize}
In the first two cases, the boom we want has length one, with a one-vertex log. In the third case, let  $u,Q,w, P'$ be as in that case. 
Since $P'$ is 
inside $P$ and hence $\phi(P')<\phi(P)$, we can apply the inductive hypothesis to $P',w$.  The result follows by concatenating $Q$ with the boom for $P',w$. This proves (1).

\bigskip

We suppose (for a contradiction) that $d_{ij}'\ne d_{ij}$ for some $i,j$ with $1\le i<j\le k$. Choose such a pair $(i,j)$ with 
$[I_i\to I_j]$ minimal. So we know:
\begin{itemize}
\item $d_{i',j'}=0$ for all $i', j'$ with $1\le i'<i$ and $i<j'<j$, and $d_{i',j'}=0$ for all $i', j'$ with $i<i'<j$ and $j<j'\le k$,
since $d$ is a demand function and $d_{ij}\ne 0$;
\item $d'_{i',j'} = d_{i',j'}$ for all $i',j'$ with $i\le i'<j'\le j$ and $(i',j')\ne (i,j)$, from the minimality of $[I_i\to I_j]$; and
\item $d'_{i',j'}=0$ for all $i', j'$ with $1\le i'\le i$ and $j\le j'\le k$ and $(i',j')\ne (i,j)$, since $d'$ is a partial $d$-demand and $d'_{i,j}<d_{i,j}$.
\end{itemize}
It follows that for $1\le i'<j'\le k$, if $\mac L_{i'j'}\ne \emptyset$, then either
\begin{itemize}
\item $i\le i'<j'\le j$; or
\item $1\le i'<j'\le i$, or $j\le i'<j'\le k$.
\end{itemize}
Let $\mac P_1, \mac P_2$ be the union of all the sets $\mac L_{i'j'}$ satisfying the first bullet and second bullet respectively.
Thus, $\mac P_1, \mac P_2$ are disjoint and have union $U(\mac L)$. For $i = 1,2$ let 
$Z_i$ be the set of all vertices with distance at most $c$ from some member of $\mac P_i$.
\\
\\
(2) {\em There is a path
from $I_i$ to $I_j$ in $G\setminus (Z_1\cup Z_2)$.}
\\
\\
Suppose that there is no such path.
Consequently
there is a region of $G\setminus (Z_1\cup Z_2)$ containing points of $[I_i\to I_j],[I_j \to I_i]$. (Note that this is true even if one of $I_i\cap V(G), I_j\cap V(G)$ is a subset of
$Z_1\cup Z_2$.) 
Hence either:
\begin{itemize}
\item there is a region of $G$ containing points $a_1\in [I_i\to I_j]$ and $a_2\in [I_j \to I_i]$; or
\item there is a region of $G$ containing a point $a_2\in [I_j \to I_i]$ and incident with a vertex $z_1\in Z_1$; or
\item there is a region of $G$ containing a point $a_1\in [I_i\to I_j]$ and incident with a vertex $z_2\in Z_2$; or
\item there is a region of $G$ incident with a vertex $z_1\in Z_1$ and with a vertex $z_2\in Z_2$.
\end{itemize}
(Note that  the points $a_1,a_2$ might lie in the interior of members of $\mac I$.)
In the first case, there is a boom of length zero joining $a_1,a_2$, and yet $(I_i, I_j)$ crosses $(a_1,a_2)$ and $d_{i,j}\ge 1$, a contradiction. 

In the second case, since $z_1\in Z_1$, there is a log $Q$ between $z_1$ and some $u \in V(P)$, for some $P\in \mac P_1$. 
Let $P\in \mac L_{i'j'}$. 
By (1),
there exists $b\in [I_{i'}\to I_{j'}]$, not in the interior of any member of $\mac I$,
such that there is a boom $Q_1\LL Q_t$
with $u\in V(Q_1)$, attaching to $b$ and with $V(Q_1)=\{u\}$, of length $t=f(V(P)\cap I_{i'},b)$. Consequently $Q, Q_2\LL Q_t$ is a 
boom of the same length $t$,
joining $a_2$ and $b$. 
Now there are two subcases. If $(i',j')\ne (i,j)$, then 
$f(V(P)\cap I_{i'},b)$ equals the sum of $d_{i''j''}$ over all choices of $i''<j''$ such that $I_{i''}\subseteq [I_i'\to b]$
and $I_{j''}\subseteq [b\to I_{j'}]$. This sum is less than the sum of $d_{i''j''}$ over all choices of $i'',j''$ with $i''<j''$ such that 
$(I_{i''},I_{j''})$ crosses $(a_2,b)$, since this latter sum also contains $d_{ij}>0$. Thus 
$t$ is less than the sum of $d_{i'',j''}$ over all $i'',j''$ such that $(I_{i''},I_{j''})$ crosses $(a_2,b)$, a contradiction. 
In the second subcase, when $(i',j')=(i,j)$, then $f(V(P)\cap I_{i'},b)$ equals $d'_{i,j}$ plus the sum of $d_{i''j''}$ over all 
choices of $i''<j''$ such that $I_{i''}\subseteq [I_i\to b]$
and $I_{j''}\subseteq [b\to I_{j}]$ and $(i'',j'')\ne (i,j)$. Again, this sum is less than 
the sum of $d_{i''j''}$ over all choices of $i''<j''$ such that 
$(I_{i''},I_{j''})$ crosses $(a_2,b)$, since this latter sum contains $d_{ij}>d'_{ij}$, and again we obtain a contradiction.

The third and fourth cases are similar, using $[I_j\to I_i]$ in place of $[I_i \to I_j]$ in the third case, and combining both arguments for the 
fourth case. This proves (2).

\bigskip

Since there is a path $P$
from $I_i$ to $I_j$ in $G\setminus (Z_1\cup Z_2)$, we can choose such a path $P$ to minimize the set of regions that belong to the side of $P$ in $\Sigma$  (in the natural sense)
that contains $[I_i\to I_j]$. 
It follows that
if we define $\mac L_{ij}'=\mac L_{ij}\cup \{P\}$, and $\mac L_{i'j'}'= \mac L_{i'j'}$ for all $(i',j')\ne (i,j)$, 
then $(\mac L_{i'j'}:1\le i'<j'\le k)$ is pushed, contrary to the maximality of the choice of $d'$.  This contradiction shows that 
$d'=d$, and so 
proves \ref{intlinkage}.~\bbox

As we said, this result is a strengthening of \ref{linkage}.
There are other strengthenings possible, that can be proved in the same way.
First, instead of specifying pairs of vertices and asking for paths
joining them with pairwise distance more than $c$, one could specify disjoint subsets of vertices, and ask for connected subgraphs,
each including one of the subsets, again with pairwise distance more than $c$; and again such subgraphs exist if and only if
there is no obstructing boom. (When $c=0$, this was done in~\cite{GM6}.)
Second, instead of asking for the paths pairwise to have distance more than $c$, one could specify a minimum distance for each pair
of paths that lie on a common region of the union of the paths, and a similar result holds (varying the lengths of logs
in a boom appropriately). Third, we could make the graph a digraph, and ask for each path to be a directed path from a specified end to its other end. We omit the details.


\section{The coarse Menger conjecture with a bounded number of alternations}

Now we turn to the main topic, the proof of \ref{mainthm}. 
The proof uses edge-contraction, which is easy enough for general graphs but difficult to be precise about when 
graphs are drawn in a surface. In all cases, edge-contraction is to be applied to the drawing, not just to the graph, and we beg the 
reader's indulgence for any imprecision. 
We are avoiding loops and parallel edges in this paper, but they might appear when we contract edges, so we need some way to get rid of them.
Let us just say that, when we do a contraction, any loops that appear are deleted, and if we make a set of parallel edges, then all but one of the edges in the set are deleted, choosing the one survivor arbitrarily. 

We are given $S,T$, subsets of vertices of a graph $G$ drawn in the plane, with $S,T$ on the outside.
With such a drawing, we can choose a simple closed curve in the plane that passes through each vertex in $S\cup T$ (and therefore 
passes through each exactly once), bounding an open 
disc in which all the remainder of the graph is drawn. Let us call such a curve a {\em bounding curve}. (We stress that the only vertices drawn in the bounding curve itself are those in $S\cup T$; all the others are strictly inside.) If $G$ has cut-vertices, or even 
worse, if $G$ is not connected, the order in which a bounding curve passes through the vertices of $S\cup T$ may not be 
unique, but that will not matter.

Given $G,S,T$ as above, if we contract some set $F$ of edges of $G$, making a graph $G'$ say, let $H$ be the subgraph of $G$ with vertex set $V(G)$ and edge set $F$. Each vertex $v$ of $G'$ is made by identifying
the vertices of some component $\eta(v)$ of $H$ under contraction, and we call $\eta(v)$ the {\em pre-image of $v$ under 
contracting $F$}. Let $S'$ be the set of $v\in V(G')$ such that 
$S\cap V(\eta(v))\ne \emptyset$, and define $T'$ similarly. Then $G'$ is also drawn in the plane and $S',T'$ are on the outside. 
We say that ``$G', S', T'$ are obtained from $G,S,T$ by contracting $F$''. 

Let $C$ be a bounding curve, and let us assign it a direction of rotation called ``clockwise''. We may assume that 
$S,T\ne \emptyset$. If $I$ is an interval of $C$ with $|I|=\{x\}$, we say $x$ is
the {\em end} of $I$ and the {\em interior} of $I$ is null. We can choose a set $\mac I=\{I_1,I_2\LL I_{2n}\}$  of 
pairwise internally-disjoint intervals of $C$
with the following properties:
\begin{itemize}
\item $I_1,I_2\LL I_{2n}$ are numbered in clockwise order on $\bd(\Sigma)$;
\item for $1\le i\le 2n$ with $i$ odd, the end(s) of $I_i$ belong to $S$, and for $i$ even, the end(s) of $I_i$ belong to $T$;
\item $I_1,I_3,I_5\LL I_{2n-1}$ are pairwise disjoint, and $I_2,I_4\LL I_{2n}$ are pairwise disjoint;
\item $S\subseteq I_1\cup I_3\cup I_5\cupcup I_{2n-1}$, and $T\subseteq I_2\cup I_4\cupcup I_{2n}$;
\end{itemize}
Thus, $\mac I$ is an interval system for $\bS$.
We call such a set $\mac I$ an {\em interval covering} for $S,T$; and we care about the size of $\mac I$. The next result shows that,
if $|\mac I|$ is bounded, we can find the $k$ connected subgraphs as in \ref{mainthm} not only with bounded diameter, but with a bounded number of edges in total.

\begin{thm}\label{boundedcover}
Let $G$ be drawn in the plane, and let $S,T\subseteq V(G)$ be on the outside. Let $C$ be a bounding curve, and let $\mac I$ be an interval covering
in $C$ of $S,T$ with size $2n$.
Suppose that there do not exist $k+1$ paths between $S$ 
and $T$, pairwise with distance more than $c$.
Then there exist $k$ connected subgraphs $B_1\LL B_{k}$, such that every $S$-$T$ 
path in $G$ contains a vertex of $B_1\cupcup B_{k}$, and $B_1\cupcup B_{k}$ has at most $8kn^2c$ edges.
\end{thm}
\Proof
Let $\mac I = \{I_1,I_2\LL I_{2n}\}$, numbered as in the definition.
Let $I_i$ have ends $a_i, b_i$, where $b_i$ is the clockwise end of $I_i$. Moreover there is an interval $L_i$
of $C$ with ends $b_i, a_{i+1}$ (where $a_{2n+1}$ means $a_1$) such that $a_{i+1}$ is the clockwise end of $L_i$.
Thus $I_1\LL I_{2n}, L_1\LL L_{2n}$ are pairwise internally-disjoint intervals of $C$ with union $C$.
For all choices of $J_1,J_2\in \{I_1\LL I_{2n}, L_1\LL L_{2n}\}$,
if there is a boom of length at most $k$ joining some point of $J_1$ and some point of $J_2$, let $F(J_1,J_2)=F(J_2,J_1)$ be the union of the 
edge-sets of the logs in some shortest such boom. If there is no such boom, let  $F(J_1,J_2)=\emptyset$. 
Thus each $F(J_1,J_2)$ has size at most $kc$. Let $F$ be the union of the sets $F(J_1,J_2)$ over all choices of $\{J_1,J_2\}$. 
Since there are at most $8n^2$ choices of the unordered pair $\{J_1,J_2\}$, 
$F$ has size at most $8kn^2c$. 

Let $G', S', T'$ be obtained from $G,S,T$ by contracting $F$. For each $v\in V(G')$, let $\eta(v)$ its pre-image under contracting $F$.
\\
\\
(1) {\em There do not exist $k+1$ vertex-disjoint paths in $G'$ from $S'$ to $T'$.}
\\
\\
Suppose that there are such paths $P_1\LL P_{k+1}$. For $1\le h\le k+1$, let $P_h$ have ends $s_h\in S$ and $t_h\in T$. 
(The choice of $s_h,t_h$ is not always uniquely determined by $P_h$, since it might have both ends in $S\cap T$, but let us choose 
some such labeling of its ends.)
For $1\le i<j\le 2n$:
\begin{itemize}
\item if $i,j$ have the same parity define $d_{ij}=0$, 
\item if $i$ is odd and $j$ is even let $d_{ij}$ be the number of $h\in \{1\LL k+1\}$ with $s_h\in I_i$ and $t_h\in I_j$; and
\item if $i$ is even and $j$ is odd, let $d_{ij}$ be the number of $h\in \{1\LL k+1\}$ with $t_h\in I_i$ and $s_h\in I_j$.
\end{itemize}
(Thus each of $P_1\LL P_{k+1}$ is counted exactly once.)
Thus $d$ is a demand function for $\mac I$. 
By hypothesis, there is no $(c+1)$-distant linkage in $G$ for $d,\mac I$.
Hence by \ref{intlinkage}, there exist $a,b\in C$, 
and a boom of length less than $t$ joining $a,b$,
where $t$ is the sum of $d_{ij}$ over $1\le i<j\le n$
such that $(I_i,I_j)$ crosses $(a,b)$ in $C$. Choose $X,Y\in\{I_1\LL I_{2n}, L_1\LL L_{2n}\}$
such that $a\in X$ and $b\in Y$. Since $t\le k+1$,
there is a boom of length at most $k$ joining $X,Y$, such that $F$ contains the edges of all its logs.
But all these edges are contracted in making $G'$; and this contradicts that $P_1\LL P_{k+1}$ exist. This proves (1).

\bigskip

From (1) and Menger's theorem, there is a set $X\subseteq V(G')$ with $|X|\le k$ such that every $S'-T'$ path in $G'$ contains a vertex of $X$. Consequently every $S$-$T$ path in $G$ contains a vertex of one of the subgraphs $\eta(x)\;(x\in X)$. But each $\eta(x)$
has all its edges in $F$. This proves \ref{boundedcover}.~\bbox

\section{The coarse Menger conjecture in the general disc case}

To complete the proof of \ref{mainthm}, we need to show how to reduce the general case to the case covered in the previous section.
If $G$ is connected and non-null and  drawn in the plane, let $r$ be its infinite region.
There is a closed walk 
$v_0,e_1,v_1\LL e_n,v_n=v_1$ that traces its boundary in the natural sense: that is, every edge incident with $r$ 
appears once or twice in this walk, and for $1\le i\le n$, the edges $e_i,e_{i+1}$ and their common end $v_i$ make an ``angle'' of 
the boundary of $r$. We call this a {\em boundary walk} of the drawing.

If $W$ is a closed walk $v_0,e_1\LL e_n, v_n = v_0$, a {\em subwalk} of $W$ is a walk of the form $v_i, e_{i+1}, v_{i+1}\LL e_j, v_j$
or of the form $v_j,e_{j+1}\LL v_n=v_0, e_1,v_1\LL v_i$, where $0\le i\le j\le k$.
In particular, if $1\le i<j\le k$, we call the first the {\em $(i,j)$-subwalk} and the second the {\em $(j,i)$-subwalk}.
We need the following lemma. 

\begin{thm}\label{repeats}
Let $G$ be a connected graph, drawn in the plane, and let $v_0,e_1,v_1\LL e_n,v_n=v_0$ be a boundary walk. Let $X\subseteq V(G)$,
such that for $1\le i<j\le n$, if $v_i=v_j$ then one of $v_{i+1}\LL v_{j-1}\in X$, and one of $v_{j+1}\LL v_n,v_1\LL v_{i-1}$ is in $X$.
Then the number of $i\in \{1\LL n\}$ with $v_i\in X$ is at most $2|X|$.
\end{thm}
\Proof
There is a hypothesis that for $1\le i<j\le n$, if $v_i=v_j$ then one of $v_{i+1}\LL v_{j-1}\in X$, and one of $v_{j+1}\LL v_n,v_1\LL v_{i-1}$ is in $X$; let us call this the {\em betweenness hypothesis}.
For any walk $W$ with terms $w_0, f_1,w_1\LL f_m, w_m=w_0$ in $G$, let $\phi(W)$ be the number of $i\in \{1\LL m\}$ such that 
$w_i\in X$. 

Now let $W$ be the given boundary walk of $G$. We proceed by induction on the length of $W$. 
If for each $x\in X$ there is at most one $i\in \{1\LL n\}$
such that $v_i=x$, then $\phi(W)\le |X|$ and the claim is true; so we assume that there exist $1\le i<j\le n$ such that 
$v_i=v_j\in X$. Choose such $i,j$
with $j-i$ minimum, and let $W_1,W_2$ be the $(i,j)$-subwalk and the $(j,i)$-subwalk respectively.
Thus there are two subdrawings $G_1,G_2$ of the drawing of $G$, such that $W_i$ is a boundary walk of $G_i$
for $i = 1,2$, and $V(G_1\cap G_2)=\{v_i\}$. Let $X_i=X\cap V(G_i)$ for $ i = 1,2$. From the minimality of $j-i$, it follows that 
$\phi(W_1)=|X_1|$; and from the hypothesis,
there exists  $h$ with $i< h< j$ such that $v_h\in X$, and so $|X_1|\ge 2$. 
Moreover, $W_2$ satisfies the betweenness hypothesis, so we can apply the inductive hypothesis to it,
and therefore  $\phi(W_2)\le 2|X_2|$. But 
$\phi(W)=\phi(W_1)+\phi(W_2)$, and $|X_1|+|X_2|=|X|+1$, and so 
$$\phi(W) \le |X_1|+ 2|X_2|=2(|X_1|+|X_2|)-|X_1|=2|X|+2-|X_1|\le 2|X|$$
since $|X_1|\ge 2$. This proves \ref{repeats}.~\bbox

We also need the following.
\begin{thm}\label{farapart}
Let $G$ be a connected graph drawn in the plane, let $C$ be a bounding curve, let $\Sigma$ be the closed disc bounded by $C$,
and let $\mac I=\{I_1\LL I_{2n}\}$ be an interval system in $C$. For 
$1\le i\le 2n$, let $P_i$ be an $S$-$T$ path with one end in $I_i$ and the other in $I_{i+1}$, such that if $P_i$ has length $>0$,
then $[I_i\to I_{i+1}]$ has nonnull interior, and all its vertices 
are incident with the region of $G$ in $\Sigma$ that contains the interior of $[I_i\to I_{i+1}]$ (where $I_{2n+1}$ means $I_1$).
Let $k\ge 1$. Then either
\begin{itemize}
\item there are $k+1$ of the paths $P_1\LL P_{2n}$ pairwise $(c+1)$-distant; or
\item there are connected subgraphs $B_1\LL B_t$ with $t\le 5k/2-1$, such that each $B_i$ has diameter at most $3c$, and each $P_i$ has a vertex in some $B_j$.
\end{itemize}
\end{thm}
\Proof We proceed by induction on $n$. ($I_{2n+1}$ means $I_1$ throughout.)
Suppose that $\dist_G(P_i,P_j)\le c$ for some non-conseutive $i,j$ with $1\le i\le j\le 2n$ (that is, with $j\ge i+2$ and $i+2n\ge j+2$). We may assume that $i = 1$.
Let $Q$ be a log between
$P_1,P_j$, and let $B$ be the maximal subgraph of $G$ such that each of its vertices has distance at most $c$ from $Q$.
Thus $B$ has diameter at most $3c$.
Let $k_1$ be the maximum number of $P_2\LL P_{j-1}$ that have distance more than $c$ from $Q$ and more than $c$ from each other.
These $k_1$ paths also have distance at least $c+1$ from each of $P_{j+1}\LL P_{2n}$, since if  some $P_{j'}$ with $j+1\le j\le 2n$
has positive length then 
all of its vertices are incident with the region of $G$ in $\Sigma$ that contains the interior of $[I_{j'}\to I_{j'+1}]$ (and so are on 
the side of $Q$ not containing $P_2$). In particular, they have distance $>c$ from $P_{2n}$, and consequently we may assume that $k_1\le k-1$.

Similarly, let $k_2$ be the maximum number of $P_{j+1}\LL P_{2n}$ that have distance more than $c$ from $Q$ and more than $c$ from
each other; and then $k_2\le k-1$. From the inductive hypothesis applied to $P_2\LL P_{j-1}$, if $k_1\ge 1$, there is a set of at most 
$5k_1/2-1$ connected subgraphs of $G$,
each of diameter at most $3c$, such that each of $P_2\LL P_{j-1}$ intersects one of them, or intersects $B$; and similarly if $k_2\ge 1$, there
is a set of at most 
$5k_2/2-1$ connected subgraphs of $G$,
each of diameter at most $3c$, such that each of $P_{j+1}\LL P_{2n}$ intersects one of them, or intersects $B$. 

Assume first that $k_1,k_2\ge 1$. 
Consequently, we have 
in total at most $5k_1/2-1 + 5k_2/2-1 +1$ connected subgraphs, each of diameter at most $3c$, such that each of $P_1\LL P_{2n}$ 
intersects one of them. If $k_1+k_2>k$ then the first outcome holds, and otherwise the second holds. 

So we may assume that $k_1=0$ say. But then all of $P_1\LL P_j$ meet $B$, and so if $k_2>0$, we have $(5k_2/2-1)+1\le 5k/2-1$
connected subgraphs that work, and if also $k_2=0$ then $B$ suffices by itself. In either case, the result holds.

So we can assume that there is no such log $Q$. Hence $P_1,P_3\LL P_{2n-1}$ are pairwise $(c+1)$-distant, and so $n\le k$. But 
we may choose $2n$ singleton sets such that each $P_i$ contains one of them, and so we may assume that $2n>5k/2-1$, and so $n=k= 1$.
If
 $\dist_G(P_1,P_2)>c$ then the first outcome holds, and otherwise there is a log joining $P_1,P_2$ and the second holds. This proves
\ref{farapart}.~\bbox

We remark that $5/2$ might not be the best possible constant in \ref{farapart}, although it is easy to show that no constant less than $3/2$ works in general.

We also need:
\begin{thm}\label{cutpoints}
Let $G$ be a connected graph drawn in the plane, with $S,T$ on the outside, and let $W$ be its boundary walk. 
Suppose that $X\subseteq V(G)$ has the property that
every subwalk of $W$ that contains a vertex in $S$ and a vertex in $T$ also contains a vertex in $X$; 
and suppose that there do not exist $k+1$ $S$-$T$ paths pairwise with distance more than $c$.
Then there exist $k$ connected subgraphs $B_1\LL B_{k}$, such that every $S$-$T$
path in $G$ contains a vertex of $B_1\cupcup B_{k}$, and $B_1\cupcup B_{k}$ has at most $32k|X|^2c$ edges.
\end{thm}
\Proof We proceed by induction on the length of $W$. Let $W$ be $v_0,e_1,v_1\LL e_n, v_n$, and suppose first that there exist 
$1\le i<j\le n$ such that $v_i=v_j$ and none of $v_{i+1}\LL v_{j-1}$ belong to $X$. It follows that at least one of $S,T$ is disjoint
from $\{v_{i+1}\LL v_{j-1}\}$. Let $W_1,W_2$ be the $(i,j)$- and $(j,i)$-subwalks respectively.
There are subdrawings $G_1,G_2$ such that $G_1\cup G_2=G$ and 
$V(G_1\cap G_2)=\{v_i\}$, such that $W_i$ is a boundary walk of $G_i$ for $i = 1,2$.
If both $S,T$ are disjoint from $\{v_{i+1}\LL v_{j-1}\}$, then deleting $V(G_1)\setminus \{v_i\}$ makes no difference, and the 
result follows by induction applied to $W_2,G_2$. We assume then that $S\cap \{v_{i+1}\LL v_{j-1}\}\ne \emptyset$, and so $T=T\cap V(G_2)$.
Let $S_2=(S\cap V(G_2))\cup \{v_i\}$. Then every subwalk of $W_2$ between $S_2,T$ contains a vertex in $X$, 
because either it is a subwalk of $W$, or it passes through $v_i$, and then it can be extended to a subwalk of $W$ by adding a portion 
of $W_1$. But then the result follows from the inductive hypothesis applied to $G_2,S_2,T,X\cap V(G_2)$. 

We may therefore assume that there is no such pair $i,j$. Consequently the hypotheses of \ref{repeats} are satisfied, and so 
the number of $i\in \{1\LL n\}$ with $v_i\in X$ is at most $2|X|$. Hence there is a partition of $\{1\LL n\}$ into at most $4|X|$
intervals (where we count a set of the form $\{j\LL n,1\LL i\}$ as an interval if $i<j$), such that for each of these intervals $I$ say,
one of $S,T$ is disjoint from $\{v_i:i\in I\}$. (To see this, take $\{i\}$ for each $v_i\in X$ as a singleton interval, and also take 
all the gaps between them.)
Let $C$ be a bounding curve. It follows (by intersecting each of our intervals with 
$V(C)\cap V(G)$) that there is an interval system $\mac I$ in $C$ of size at most $4|X|$, such that 
$S\cup T\subseteq \bigcup_{I\in \mac I} I$, and such that for each $I\in \mac I$, one of $S,T$
is disjoint from $I$. By replacing by their union any two 
consecutive intervals in $\mac I$ that are both disjoint from $S$ or both disjoint from $T$, we may assume that $\mac I$ is an 
interval covering of $S,T$ of size at most $4|X|$. Hence the result follows from \ref{boundedcover}. This proves \ref{cutpoints}. ~\bbox

Now we deduce our main theorem, which we restate:
\begin{thm}\label{mainthm2}
Let $c,k\ge 0$, let $G$ be drawn in the plane, with $S,T$ on the outside.
Then either:
\begin{itemize}
\item there are $k+1$ paths between $S,T$, pairwise at distance more than $c$; or
\item there is a set of at most $k$ connected subgraphs of $G$, 
such that every path between $S,T$ intersects one of these subgraphs, and the sum of the diameters of the subgraphs is at most 
$200k^3c$.
\end{itemize}
\end{thm}
\Proof
We may assume that $G$ is connected. Let $C$ be a bounding curve, let $\Sigma$ be the closed disc bounded by $C$,
and let $\mac I=\{I_1\LL I_{2n}\}$ be an interval covering of $S,T$. For 
$1\le i\le 2n$, let $P_i$ be an $S$-$T$ path with one end in $I_i$ and the other in $I_{i+1}$, such that for each $v\in V(P_i)$,
$v$ is incident with a region of $G$ in $\Sigma$ that contains a  point of $[I_i\to I_{i+1}]$. By \ref{farapart}, 
we may assume that there is a set of at most $5k/2-1$ 
connected subgraphs of $G$,
each of diameter at most $3c$, such that each of $P_1\LL P_{2n}$ intersects one of them. These subgraphs need not be vertex-disjoint,
but their union has at most $5k/2$ components, and the sum of the diameters of these components is at most $15kc/2$. 
Let $F$ be the union of the 
edge sets of these components, and let $G',S',T'$ be obtained from $G,S,T$ by contracting $F$. Each of the components with union $F$
contracts to a vertex in $G'$; let $X$ be the set of such vertices. Thus $|X|\le 5k/2-1$. For each $v\in V(G')$,
let $\eta(v)$ be the pre-image of $v$ under contracting $F$. Thus the sum of the diameters of $\eta(v)$, over all $v\in V(G')$,
is at most $15kc/2$. (Note that, typically, $\eta(v)$ consists just of the vertex $v$ and so has diameter zero: the only terms that contribute to the sum are when $v\in X$.)

Now $G'$ is connected; let $W$ be a boundary walk, where $W$ is $v_0,e_1\LL e_t,v_t=v_0$ say. Every subwalk of $W$ that contains a 
vertex of $S'$ and a vertex of $T'$,
also contains a vertex of $X$, because of the way we constructed $X$. We may assume that there do not exist 
$k+1$ paths of $G'$ between $S',T'$, pairwise at distance more than $c$, since otherwise the first outcome of the theorem holds. 
Thus, by \ref{cutpoints}, 
there exist $k$ connected subgraphs $B_1'\LL B_{k}'$ of $G'$, such that every $S'-T'$
path in $G'$ contains a vertex of $B_1'\cupcup B_{k}'$, and $B_1'\cupcup B_{k}'$ has at most $32k(5k/2-1)^2c$ edges.
For $1\le i\le k$, let $B_i$ be the subgraph of $G$ formed by the union of the edges of $B_i'$ and the subgraphs 
$\eta(v)\;(v\in V(B_i'))$. Thus $B_1\LL B_{k}$ are connected,
and every $S$-$T$ path in $G$ meets one of them. Moreover, the sum of the diameters of $B_1\LL B_{k}$ is at most 
$32k(5k/2-1)^2c+ 15kc/2\le 200k^3c$, since $B_1'\cupcup B_{k}'$  has at most $32k(5k/2-1)^2c$ edges and the sum of the diameters of $\eta(v)$ for $v\in V(G')$
is at most $15kc/2$. This proves \ref{mainthm2}.~\bbox

\section{An algorithm}

Suppose we are given a graph drawn in the plane, with $S,T$ on the outside. Can we check in polynomial time whether there exist $k$
$S$-$T$ paths pairwise $(c+1)$-distant? If $c,k$ are constants, the answer is yes, as follows. Let us say the {\em depth} of a vertex $v$
is the minimum $n$ such that there is a sequence 
$$v=v_0,r_1,v_1,r_2\LL v_n$$
of alternating vertices and regions, where $v_n$ is incident with the infinite region, and $r_i$ is incident with $v_{i-1}, v_i$ 
for $1\le i\le n$.
The $S$-$T$ paths we want exist if and only if there exist distinct $s_1\LL s_k$ in $S$ and distinct $t_1\LL t_k\in T$,
such that there are $k$ paths joining $s_i,t_i$ for $1\le i\le k$, pairwise $c$-distant. And for a given choice of 
$s_1\LL s_k, t_1\LL t_k$, whether such paths exist is determined by the existence of certain booms of length at most $k-1$, 
attaching to two points of some bounding curve, as in \ref{linkage}. The vertices in every such boom have depth 
at most $c(k-1)/2$, so we can delete all vertices with depth more than $c(k-1)/2 +1$ without changing whether the booms exist, and therefore without changing whether the paths we want exist. 
(The ``$+1$'' is to avoid making new regions incident with vertices at depth $c(k-1)/2$.)
But after this, the graph has bounded tree-width, and the problem can be solved in linear time, by Courcelle's theorem~\cite{courcelle},
since the question can be expressed by a monadic second-order formula.

\end{document}